\documentclass[11pt,a4paper]{article}

\usepackage{times}
\usepackage[utf8]{inputenc}
\usepackage{latexsym,amssymb,amsmath}
\usepackage[pdftex,colorlinks=true,breaklinks=true]{hyperref}
\usepackage{xspace}
\usepackage{tikz}
\usepackage[pdftex,colorlinks=true]{hyperref}

\newcommand{\Galg}{\mathbf{G}}
\newcommand{\Talg}{\mathbf{T}}
\newcommand{\F}{\mathbb{F}}
\newcommand{\Fp}{\mathbb{F}_p}
\newcommand{\overFp}{\overline{\mathbb{F}}_p}

\newcommand{\overFq}{\overline{\mathbb{F}}_q}
\newcommand{\Qpp}{\mathbb{Q}_{p'}}
\newcommand{\Z}{\mathbb{Z}}
\newcommand{\C}{\mathbb{C}}
\newcommand{\Hom}{\operatorname{Hom}}
\newcommand{\Spin}{\operatorname{Spin}}
\newcommand{\Sp}{\operatorname{Sp}}
\newcommand{\SL}{\operatorname{SL}}
\newcommand{\CHEVIE}{\textsf{CHEVIE}\xspace}
\newcommand{\ATLAS}{\textsf{ATLAS}\xspace}
\newcommand{\CTblLib}{\textsf{CTblLib}\xspace}
\newcommand{\AtlasRep}{\textsf{AtlasRep}\xspace}
\newcommand{\MeatAxe}{\textsf{MeatAxe}\xspace}
\newcommand{\Magma}{\textsf{Magma}\xspace}
\newcommand{\GAP}{\textsf{GAP}\xspace}
\newcommand{\id}{\textrm{id}}

\newtheorem{Thm}{Theorem}[subsection]
\newcommand{\eprf}{\hfill$\Box$}

\begin{document}

\author{Frank Lübeck}
\title{Turning  Weight  Multiplicities  into  Brauer  Characters}
\date{August 2018}
\maketitle

\begin{abstract}
We describe  two methods for  computing $p$-modular Brauer  character tables
for groups  of Lie  type $G(p^f)$ in  defining characteristic  $p$, assuming
that  the ordinary  character table  of $G(p^f)$  is known,  and the  weight
multiplicities of  the corresponding algebraic  group $\Galg$ are  known for
$p$-restricted highest weights.
\end{abstract}

\section{Introduction}\label{sec-intro}

Let $G(q)$  be finite group of  Lie type arising from  a connected reductive
algebraic group $\Galg$  of simply-connected type over  an algebraic closure
$\overFq$ of characteristic $p$.

The   irreducible  representations   of   $G(q)$   over  $\overFp$   (called
defining  characteristic representations)  are  restrictions of  irreducible
representations   of   the   algebraic  group   $\Galg$.   The   irreducible
representations  of $\Galg$  over $\overFp$  can  be described  in terms  of
weight multiplicities.

If   the   weight   multiplicities   for  $\Galg$   are   explicitly   known
for   sufficiently  many   irreducible  representations   (parameterized  by
$p$-restricted highest weights)  a lot of information  about the irreducible
representations of $G(q)$  for small powers $q$ of $p$  can be computed. For
example  the degrees  of all  irreducible representations  of $G(q)$  or the
composition factors of tensor products of irreducible representations.

But  this  description  does  not   yield  a  relation  of  the  irreducible
representations of $G(q)$ in characteristic  $p$ with the ordinary (complex)
representations  of $G(q)$.  This relation  is provided  by the  $p$-modular
Brauer  character   table  of   $G(q)$,  equivalently  by   the  $p$-modular
decomposition  matrix  of  $G(q)$  or   the  (ordinary)  characters  of  the
$p$-modular projective indecomposable modules.

In this note we describe two ways to compute from weight multiplicities of a
representation of $\Galg$ the corresponding Brauer character of $G(q)$.

As an  application of  these methods  we compute  some character  tables for
the  modular {\ATLAS}  project~\cite{ModAtl,WebModAtl},  which  has the  aim
to  provide all  Brauer  character  tables for  all  primes  for the  groups
whose ordinary  character table  is mentioned in  the {\ATLAS}~\cite{ATLAS}.
More   precisely,  we   compute  the   full  $p$-modular   Brauer  character
tables  for  the groups  (in  {\ATLAS}  notation) $F_4(2)$,  $2^2.O_8^+(3)$,
$2.O_8^-(3)$, $O_{10}^+(2)$, $O_{10}^-(2)$ and  $S_{10}(2)$ where $p$ is the
defining  characteristic.  We also  find  partial  tables for  $E_6(2)$  and
$3.{}^2\!E_6(2)$.

One possible  application of  these tables  is the  induction of  Brauer and
projective characters to (sporadic) overgroups.

In Section~\ref{sec-setup} we  describe the groups we consider  here in more
detail and explain how we represent them for computations.

In  Section~\ref{sec-mults} we  give  an overview  of weight  multiplicities
encoded   in   dominant   characters,    and   in   particular   we   sketch
in~\ref{ss_computechars} how to compute  with such dominant characters (this
may be of independent interest).

Section~\ref{sec-smallrep}  describes  a  first  method  to  compute  Brauer
characters which  uses explicit small  degree representations of  $G(q)$ and
computations with dominant characters.

In  Section~\ref{sec-sscl} we  show  how to  compute  a parameterization  of
conjugacy classes of  elements of $p'$-order in $G(q)$  (without an explicit
representation of the group).

Section~\ref{sec-brauer}  describes  a  second   method  to  compute  Brauer
characters which  uses the  weight multiplicities directly  as well  as some
invariants of conjugacy classes.

In all sections we illustrate  the theoretical descriptions with the example
of $\Galg$ of type $D_4$ and $G(q) \in \{\Spin_8^+(3), \Spin_8^-(3)\}$.

\textbf{Acknowledgements.} I thank Klaus Lux for convincing me that for some
applications is  it not  enough to know  the irreducible  representations of
$G(q)$ in  defining characteristic only  in terms of  weight multiplicities.
This motivated the preparation of these notes. And I thank Thomas Breuer for
answering some questions about Brauer character tables, for helping with the
computation  of  some  ordinary  character tables  with  {\Magma},  and  for
distributing the new character tables described here with the {\GAP}-package
{\CTblLib}.

\section{Setup and Notation}\label{sec-setup}

\subsection{Root Data}\label{ss-rootdata}

We consider  connected reductive  groups $\Galg$  over an  algebraic closure
$\overFp$ of a finite field of characteristic $p$ and their finite subgroups
$G(q) =  \Galg^F$ of fixed points  under a Frobenius morphism  $F$.

If  $\Galg$ is  of  rank $r$  and  $\Talg  \leq \Galg$  is  a maximal  torus
of  $\Galg$  there  is  an  associated  root  datum  $(X,\Phi,Y,\Phi^\vee)$,
see~\cite[7.4, 9.6]{Sp98}. Here  $X  \cong  \Z^r$  is  the  character  group
$\Hom(\Talg, \overFp^\times)$  and $Y \cong  \Z^r$ is the  cocharacter group
$\Hom(\overFp^\times, \Talg)$ of the torus  $\Talg$, and $\Phi \subset X$ is
the set of roots and $\Phi^\vee \subset  Y$ is the set of coroots of $\Galg$
with respect to $\Talg$.

If $F$ is a Frobenius morphism of $\Galg$, we assume that $F(\Talg) = \Talg$
and in that case $F$ naturally induces $\Z$-linear maps on $X$ and $Y$ which
are of the form  $q \cdot F_0$ where $F_0$ is of  finite order permuting the
roots and $q$ is a power of $p$, see~\cite[1.17--1.18]{Ca85}. We write $G(q)
=  \Galg^F$  for the  finite  group  of  $F$-fixed  points to  indicate  the
parameter $q$ determined by $F$. (We  exclude here more general $F$, leading
for  example to  Suzuki and  Ree groups,  to keep  some statements  simpler;
except for some remarks in the proof of~\ref{ThmMults}.)

The algebraic  group $\Galg$  is determined  up to  isomorphism by  the root
datum and the  field $\overFp$. The finite group $G(q)$  of $F$-fixed points
is determined up to isomorphism by the root datum, $F_0$ and the number $q$.

For computations  we use the  setup from the  \CHEVIE~\cite{CHEVIE} software
package. See also~\cite[Section 2]{BL13} for more details.

Let $l$ be the  semisimple rank of $\Galg$. We encode the  root datum by two
matrices  $(A,A^\vee)$ in  $\Z^{l \times  r}$.  There is  a natural  pairing
$\langle  \cdot, \cdot  \rangle:  X \times  Y  \to \Z$  and  we choose  dual
$\Z$-bases of $X$  and $Y$ with respect  to that pairing. We write  a set of
simple roots with respect to the basis  of $X$ and put the coefficients into
the rows  of $A$, similarly we  write the corresponding simple  coroots with
respect  to the  basis of  $Y$ and  put the  coefficients into  the rows  of
$A^\vee$.

If $\alpha_i$ is a simple  root and $\alpha_i^\vee$ the corresponding simple
coroot ($1 \leq i \leq l$), we define
\[ s_i: X \to X, x \mapsto x - \langle x, \alpha^\vee \rangle \alpha, \]
\[ s_i^\vee: Y \to Y, y \mapsto y - \langle \alpha, y \rangle \alpha^\vee.\]
The $s_i$, $1  \leq i \leq l$, are  a set of Coxeter generators  of the Weyl
group $W$  of $\Galg$ (and the  same is true  for the $s_i^\vee$, $1  \leq i
\leq l$). Using the $i$-th rows of  $A$ and $A^\vee$ we can write down $s_i$
as matrix acting on $X$ with respect  to our basis of $X$. The transposed of
that matrix represents the action of $s_i^\vee$ on $Y$.

The set of all  roots $\Phi$ can be found as the orbits  of the simple roots
under the action of $W$.

If $F_0: X \to  X$ is a $\Z$-linear map of finite order  which fixes the set
of roots and  whose dual map $F_0^\vee:  Y \to Y$ fixes the  set of coroots,
then there is for each power $q$ of $p$ a Frobenius map $F$ on $\Galg$ which
induces $q F_0$ on $X$.

\subsection{Simply-Connected Groups}\label{ss-sc}

From now  we assume that $\Galg$  is a semisimple group  of simply connected
type. In  our setup  this means  that $r=l$  and that  the coroots  span the
lattice $Y$.  Therefore, we choose  a set of simple  coroots $\alpha_1^\vee,
\ldots, \alpha_l^\vee$  as $\Z$-basis of  $Y$. The corresponding  dual basis
$\omega_1,  \ldots,  \omega_l$ of  $X$  is  called  the set  of  fundamental
weights. The  matrix $A^\vee$  is the  identity matrix of  size $l$  and the
matrix $A$  is the  transposed of the  Cartan matrix of  the root  system of
$\Galg$  (the $l  \times  l$ matrix  with  $(i,j)$-entry $\langle  \alpha_j,
\alpha_i^\vee \rangle$).

The set $X_+ = \{\sum_i a_i \omega_i\mid\; a_i \in \Z_{\geq 0}\textrm{ for }
1 \leq  i \leq l\}$ is  called  the  set of dominant weights and  for a real
number $b$  we call  $X_b =  \{ \sum_i a_i  \omega_i \in  X_+\mid\; a_i  < b
\textrm{ for  } 1 \leq  i \leq l\}$ the  set of $b$-restricted  weights. The
$W$-orbit of each $\lambda \in X$ contains a unique dominant weight.

There is  a partial  order on  $X$ defined  by: $\lambda,  \mu \in  X$, then
$\lambda \geq \mu$ if and only if $\lambda - \mu = \sum_i b_i \alpha_i$ with
$b_i \in \Z_{\geq 0}$ for $1 \leq i \leq l$.

\subsection{Example}\label{ex-d4}

We consider the root system of type $D_4$. Its Cartan matrix is
\[{\scriptsize \left(\begin{array}{rrrr}%
2&0&-1&0\\%
0&2&-1&0\\%
-1&-1&2&-1\\%
0&0&-1&2\\%
\end{array}\right)} \]
and can be encoded in the Dynkin diagram
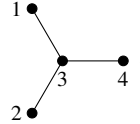
\begin{tikzpicture}[scale=0.8, baseline=-1.0]
\coordinate (3) at (0,0);
\coordinate (1) at (120:1);
\coordinate (2) at (240:1);
\coordinate (4) at (1,0);
\foreach \p in {1, 2, 3, 4} \filldraw [black] (\p) circle [radius=0.08];
\draw (1) -- (3);
\draw (2) -- (3);
\draw (3) -- (4);
\draw (1) node[anchor=east] {{\scriptsize 1}};
\draw (2) node[anchor=east] {{\scriptsize 2}};
\draw (3) node[anchor=north] {{\scriptsize 3}};
\draw (4) node[anchor=north] {{\scriptsize 4}};
\end{tikzpicture}.

We encode  the root datum of  the simply connected algebraic  groups of type
$D_4$ by the pair of matrices
\[ (A, A^\vee) = ({\scriptsize
\left(\begin{array}{rrrr}%
2&0&-1&0\\%
0&2&-1&0\\%
-1&-1&2&-1\\%
0&0&-1&2\\%
\end{array}\right)},
{\scriptsize
\left(\begin{array}{rrrr}%
1&0&0&0\\%
0&1&0&0\\%
0&0&1&0\\%
0&0&0&1\\%
\end{array}\right) } ).\]
Then the generating reflections $s_1, \ldots, s_4$ of the Weyl group $W$ act
on $X$ by
\[{\scriptsize
\left(\begin{array}{rrrr}%
-1&0&1&0\\%
0&1&0&0\\%
0&0&1&0\\%
0&0&0&1\\%
\end{array}\right)%
,
\left(\begin{array}{rrrr}%
1&0&0&0\\%
0&-1&1&0\\%
0&0&1&0\\%
0&0&0&1\\%
\end{array}\right)%
,
\left(\begin{array}{rrrr}%
1&0&0&0\\%
0&1&0&0\\%
1&1&-1&1\\%
0&0&0&1\\%
\end{array}\right)%
,
\left(\begin{array}{rrrr}%
1&0&0&0\\%
0&1&0&0\\%
0&0&1&0\\%
0&0&1&-1\\%
\end{array}\right)%
}\]
and $s_1^\vee, \ldots, s_4^\vee$ act on $Y$ by the transposed of the printed
matrices.

For  each   characteristic  $p$  these  data   determine  the  corresponding
reductive  algebraic group  $\Galg$ which  is isomorphic  to the  Spin group
$\Spin_8(\overFp)$.

If we choose  as $F_0$ the identity  matrix, then for each power  $q$ of $p$
the  corresponding Frobenius  morphism of  $\Galg$ yields  a group  of fixed
points $G(q) = \Spin_8^+(q)$.

The $F_0$ that permutes the first  two simple roots (and coroots) yields the
finite groups  $G(q) =  \Spin_8^-(q)$, and from  $F_0$ permuting  the simple
roots  $1 ,2,  4$ cyclically  we find  the Steinberg  triality groups  $G(q)
={}^3\!D_4(q)$.

\section{Weight Multiplicities}\label{sec-mults}

We use the setup from~\ref{ss-sc}. So $\Galg$ is a connected reductive group
of  simply-connected type  over  $\overFp$ with  a  Frobenius morphism  $F$,
maximal torus $\Talg  = F(\Talg)$, and finite group of  fixed points $G(q) =
\Galg^F$.

\subsection{Irreducible Representations}\label{ss-irr}

We recall a few standard  facts about irreducible representations of $\Galg$
and $G(q)$ in their defining characteristic $p$.

If  $M$   is  a  rational  module   for  $\Galg$  then  $M$   considered  as
$\Talg$-module splits into a direct  sum of non-zero subspaces $M_\mu$ which
are  common eigenspaces  for all  $t \in  \Talg$. Here  $\mu \in  X$ is  the
homomorphism $\mu:\Talg \to \overFp^\times$ which maps each $t \in \Talg$ to
its eigenvalue  $\mu(t)$ on $M_\mu$. The  set $\{\mu \in X\mid\;  M_\mu \neq
0\}$ is called  the weights of $M$, it  is a union of $W$-orbits  on $X$ and
when $\mu, \nu  \in X$ are in  the same $W$-orbit, then  $M_\mu$ and $M_\nu$
have the same dimension.

If  $M$ is  irreducible then  it has  a unique  highest weight  $\lambda \in
X_+$  with $\lambda  \geq \mu$  for all  weights of  $M$. Chevalley~\cite[II
2.7]{Ja03} showed  that such an $M$  is determined up to  isomorphism by its
highest weight,  we denote this  module $L(\lambda)$. And for  each dominant
$\lambda \in X_+$ there is such an irreducible module $L(\lambda)$.

Steinberg~\cite[13.1]{St68}  showed that  for fixed  characteristic $p$  all
$L(\lambda)$ can be described in terms of  a finite number of them. For this
write an arbitrary dominant $\lambda \in X_+$ as a finite linear combination
$\lambda = \sum_{i  = 0}^k p^i \lambda_i$ where all  $\lambda_i \in X_p$ are
$p$-restricted (write  the entries  of $\lambda$  in base  $p$). Steinberg's
tensor product theorem says that
\[ L(\lambda) = L(\lambda_0) \otimes L(\lambda_1)^{[1]} \otimes \cdots
\otimes L(\lambda_k)^{[k]},\]
where  $L(\lambda)^{[i]}$ is  the  $i$-th Frobenius  twist of  $L(\lambda)$,
$\Galg$ acts on this module by the composition of its action on $L(\lambda)$
and the field automorphism $c \mapsto c^{p^i}$ of $\overFp$.

Steinberg~\cite[13.3]{St68}  also  showed  that   the  restrictions  of  all
$\{L(\lambda)\mid\; \lambda \in X_q\}$  (the $L(\lambda)$ for $q$-restricted
weights) to $G(q)$ yield a set of representatives of all isomorphism classes
of irreducible representations of $G(q)$ over $\overFp$.

\subsection{Characters by Weight Multiplicities}\label{ss-mults}

A lot of information about  a finite dimensional rational $\Galg$-module $M$
is encoded  in the list of  weights $\mu$ of  $M$ and the dimensions  of the
weight spaces $M_\mu$. Using that the  set of weights consists of $W$-orbits
which each contain a unique dominant weight, we get an efficient description
of these data by a dominant character, which only lists the dominant weights
of $M$ and the dimensions of the corresponding weight spaces. If $\mu \in X$
is a weight of $M$ we call  $m_\mu(M) = \dim(M_\mu)$ the multiplicity of the
weight $\mu$ in $M$.

The orbit  lengths of  weights are  not difficult to  compute because  for a
weight $\sum_{i=1}^l  a_i \omega_i$ its  stabilizer in $W$ is  the parabolic
subgroup $\langle s_i \mid\; a_i = 0 \rangle$.

We  can find  all weights  in a  $W$-orbit by  a standard  orbit calculation
using  the matrices  for  the generators  $s_1, \ldots,  s_l$  of $W$  given
in~\ref{ss-rootdata}.

Here is an  example. Let $\Galg$ be of  type $D_4$, $p = 3$,  and $\lambda =
\omega_2+2 \omega_4$. We denote elements of $X$ by their coefficient vectors
with  respect  to the  basis  $\omega_1,  \ldots, \omega_4$  of  fundamental
weights  as in  Example~\ref{ex-d4}.  In  the following  table  we show  the
dominant character of $L(\lambda)$ and the  lengths of the $W$-orbits of the
mentioned dominant weights.
\[
\begin{array}{rrr}
\mu & m_\mu & |\mu^W| \\
\hline
( 0, 1, 0, 2 ) & 1 & 32\\
( 0, 1, 1, 0 ) & 1 & 48\\
( 1, 0, 0, 1 ) & 3 & 32\\
( 0, 1, 0, 0 ) & 6 & 8
\end{array}
\]

In our  example, we can use  the $s_i$ from~\ref{ex-d4} (matrices  acting on
$X$), the orbit of $(0,1,0,0)$ consists of the  weights $( 0, 1, 0, 0 )$, $(
0, -1, 1, 0 )$, $( 1, 0, -1, 1 )$, $(  -1, 0, 0, 1 )$, $( 1, 0, 0, -1 )$, $(
-1, 0, 1, -1 )$, $( 0, 1, -1, 0  )$, $( 0, -1, 0, 0 )$ and the corresponding
weight spaces of $L( (0,1,0,2) )$ have dimension $6$.

\subsection{Computing with Dominant Characters}\label{ss_computechars}

Given dominant characters of rational $\Galg$-modules $M$ and $M'$ we can

\begin{itemize}
\item[(a)]  compute  the total  dimension,  $\dim  M  = \sum_{\mu  \in  X_+}
m_\mu(M) \cdot |\mu^W|$,
\item[(b)]  compute the  dominant character  of the  $i$-th Frobenius  twist
$M^{[i]}$, for this just multiply all weights by $p^i$,
\item[(c)] compute  the character of $M  \otimes M'$, using that  for $v \in
M_\mu$  and  $v'  \in  M'_{\mu'}$  the  vector $v  \otimes  v'$  is  in  the
weight  space $(M  \otimes  M')_{\mu  + \mu'}$  because   $\mu(t) \mu'(t)  =
(\mu+\mu')(t)$ for  $t \in \Talg$;  to compute the weight  multiplicities of
$M_\mu \otimes M'_{\mu'}$ one needs to  enumerate at least one of the orbits
$\mu^W$ or $\mu'^W$,
\item[(d)]  find  the  composition  factors of  $M$  provided  the  dominant
characters of those  are available, for this enumerate  the dominant weights
of $M$ with  respect to a total ordering which  refines the partial ordering
$>$ on $X$,  then the first weight $\lambda$ and  its multiplicity $m$ shows
the multiplicity of $L(\lambda)$ as  composition factor of $M$, subtract $m$
times the dominant character of $L(\lambda)$ and proceed recursively.
\end{itemize}

Now  assume  that  we know  for  a  fixed  $p$  the dominant  characters  of
$L(\lambda)$ for all  $p$-restricted $\lambda \in X_p$. Fix the  map $F_0: X
\to X$ and let $q$ be a power of $p$.

Then we can

\begin{itemize}
\item[(e)]  compute the  dimensions  of all  irreducible representations  of
$G(q)$  over  $\overFp$, using  the  Steinberg  tensor product  theorem  for
all  $q$-restricted weights  (only the  dimensions of  the $L(\lambda)$  for
$p$-restricted $\lambda$ are needed),
\item[(f)]  compute   the  dominant  characters  of   $L(\lambda)$  for  all
$q$-restricted weights, using the tensor product theorem and~(b) and~(c),
\item[(g)] find the  composition factors of the  restriction of $L(\lambda)$
to $G(q)$ for arbitrary $\lambda \in X$, for this note that $L(\lambda)$ and
$L(\lambda  (q  F_0))$  restrict  to  the  same  representation  of  $G(q)$,
so  decompose  an  arbitrary  $\lambda  \in  X$  as  $\sum_{i  \geq  0}  q^i
\lambda_i F_0^i$  with all $\lambda_i  \in X_q$ being $q$-restricted  to get
$L(\lambda)\mid_{G(q)}  =  \bigotimes_{i \geq  0}  L(\lambda_i)\mid_{G(q)}$,
use~(h) and this point~(g) recursively,
\item[(h)]  decompose  tensor  products of  irreducible  representations  of
$G(q)$ over $\overFp$, use~(c),~(d) and~(g) for factors not corresponding to
$q$-restricted weights.
\end{itemize}

We continue the example of $\Galg$ in type $D_4$, $p=3$ from~\ref{ss-mults}.
Let $\lambda = \omega_2+2\omega_4$ and  $\lambda' = 2\omega_4$. The dominant
character of $L(\lambda)\otimes L(\lambda')$ is

\[
\begin{array}{rrr}
\mu & m_\mu & |\mu^W| \\
\hline
( 0, 1, 0, 4 ) & 1 &   32 \\
( 0, 1, 1, 2 ) & 2 &   96 \\
( 0, 1, 2, 0 ) & 3 &   48 \\
( 1, 2, 0, 1 ) & 4 &   96 \\
\end{array}
\quad
\begin{array}{rrr}
\multicolumn{3}{c}{\textit{(cont.)}}\\
\hline
( 1, 0, 0, 3 ) & 6 &   32 \\
( 1, 0, 1, 1 ) & 11 &   96 \\
( 2, 1, 0, 0 ) & 15 &   32 \\
( 0, 3, 0, 0 ) & 6 &   8 \\
\end{array}
\quad
\begin{array}{rrr}
\multicolumn{3}{c}{\textit{(cont.)}}\\
\hline
( 0, 1, 0, 2 ) & 24 &   32 \\
( 0, 1, 1, 0 ) & 34 &   48 \\
( 1, 0, 0, 1 ) & 63 &   32 \\
( 0, 1, 0, 0 ) & 112 &   8
\end{array}
\]
So,  for  the  algebraic  group  $\Galg$ the  tensor  product  contains  one
composition  factor  isomorphic  to  $L((0,1,0,4))  =  L((0,1,0,1))  \otimes
L((0,0,0,1))^{[3]}$.  This composition  factor restricted  to $\Spin_8^+(3)$
or  to  $\Spin_8^-(3)$  is  the   tensor  product  of  the  restrictions  of
$L((0,1,0,1))$ and $L((0,0,0,1))$ to those  groups, while the restriction to
$^3\!D_4(3)$ is the tensor product of the restrictions of $L((0,1,0,1))$ and
$L((0,1,0,0))$ (because  $F_0$ permutes $\omega_1 \to  \omega_2 \to \omega_4
\to \omega_1$) to  the finite group. A similar computation  shows that these
tensor  products  for  the  finite  groups  have  five  composition  factors
(corresponding to $3$-restricted weights) in each case.

\subsection{Some Dominant Characters}

\begin{Thm}\label{ThmMults}
Let $p$ and $\Galg$ as in one row of the following table
\[ \begin{array}{lll}
p& \textrm{Lie type of } \Galg& \textrm{group name}\\
\hline
3 & D_4 & \Spin_8(\overFp)\\
2 & B_4 & \Spin_9(\overFp)\\
2 & F_4 & F_4(\overFp)\\
2,3 & A_5 & \SL_{6}(\overFp)\\
2 & C_5 & \Sp_{10}(\overFp)\\
2 & D_5 & \Spin_{10}(\overFp)\\
\end{array}
\]

Then  the dominant  characters of  the irreducible  rational representations
$L(\lambda)$ of  $\Galg$ are  known for  all $p$-restricted  highest weights
$\lambda$.

If $\Galg$  is simply  connected of type  $E_6$ and $p=2$  then we  know the
dominant  characters  of  44 irreducible  $L(\lambda)$  with  $2$-restricted
$\lambda$ (there are 64 such $\lambda$).

\end{Thm}

\textrm{Proof.}
These   characters   were   computed   with  the   strategy   and   programs
described  in~\cite{luebeckdef}. The  characters  are available  on the  web
page~\cite{WebMults}. (Because of the size of these data we do not reproduce
them within  this article.)  The computations involved several weeks  of CPU
time.

We remark that the  case $F_4$ and $p=2$ is particularly  easy. In this case
$\Galg$  has an  exceptional automorphism  $\tilde  F$ whose  square is  the
Frobenius morphism which acts as $2 \; \textrm{Id}$ on $X$.

If we number the simple roots and the fundamental weights  according to  the
Dynkin diagram
\begin{tikzpicture}[scale=0.8, baseline=-3.0]
\coordinate (1) at (0,0);
\coordinate (2) at (1,0);
\coordinate (3) at (2,0);
\coordinate (4) at (3,0);
\coordinate (5) at (1.3,0.3);
\coordinate (6) at (1.3,-0.3);
\coordinate (7) at (1.7,0);

\draw[line width=0.6] (1) -- (2);
\draw[double, double distance between line centers=3.38, line width=0.6] (2) -- (3);
\draw[line width=0.6] (3) -- (4);
\draw[line width=0.6] (5) -- (7);
\draw[line width=0.6] (6) -- (7);
\draw (1) node[anchor=north] {{\scriptsize 1}};
\draw (2) node[anchor=north] {{\scriptsize 2}};
\draw (3) node[anchor=north] {{\scriptsize 3}};
\draw (4) node[anchor=north] {{\scriptsize 4}};
\foreach \p in {1, 2, 3, 4} \filldraw [black] (\p) circle [radius=0.08];
\filldraw [black] (7) circle [radius=0.001];
\end{tikzpicture}, 
then $\tilde F$ maps $\omega_1,\omega_2,\omega_3,\omega_4$ to $2 \omega_4, 2
\omega_3,\omega_2,\omega_1$,  respectively.  The  Steinberg  tensor  product
theorem also holds for $\tilde F$: Let
\[M = \{(0,0,0,0),(0,0,0,1),(0,0,1,0),
(0,0,1,1)\}\]
and write an arbitrary dominant $\lambda \in X^+$ as $\lambda = \sum_{i=0}^k
\lambda_i {\tilde F}^i$ where all $\lambda_i \in M$. Then
\[L(\lambda)  =  L(\lambda_0)   \otimes  L(\lambda_i)^{[1]}  \otimes  \cdots
\otimes L(\lambda_k)^{[k]}\]
where   now  ${}^{[i]}$   denotes  the   $i$-th  twist   with  $\tilde   F$,
see~\cite[11.2]{St68}.  It is  easy to  compute the  dominant characters  of
$L(\lambda)$ for $\lambda \in M$. This was already done in~\cite{V70}.
\mbox{}\eprf

\section{Brauer Table from Small Representation and
Tensoring}\label{sec-smallrep}

Let $\Galg$, $p$ and  $G(q)$ be as in previous sections.  In this section we
describe a first method to find the Brauer character table of a finite group
$G(q)$ in defining characteristic $p$.

We  recall that  Brauer  characters  are defined  relative  to an  embedding
of  the multiplicative  groups  $\overFp^\times \hookrightarrow  \C^\times$.
There  is  a  convention how  to  choose  this  map  which is  used  in  the
Modular  Atlas~\cite{ModAtl}  and in  the  {\CTblLib}  library of  character
tables~\cite{CTblLib}  distributed  with  {\GAP}~\cite{GAP4} (based  on  the
notion of Conway  polynomials). To compute the Brauer character  value for a
group element  given by its representing  matrix over a finite  extension of
$\Fp$ one computes its eigenvalues in  $\overFp$, lifts them to $\C$ via the
mentioned embedding, and sums up the  images. {\GAP} provides a function for
this computation. Brauer  character values are only computed  on elements of
$p'$-order.

To   start,  we   assume  that   we   have  an   explicit  faithful   matrix
representation of  $G(q)$ over a finite  extension of $\Fp$. For  the groups
considered  in  this  article  these  are  available  from  {\GAP}  commands
like  \texttt{SL(6,3)}  or \texttt{Sp(10,2)}  or  from  the ATLAS  of  group
representations~\cite{WWWAtlas},  accessible in  {\GAP} via  the {\AtlasRep}
package~\cite{AtlasRep}.

We also assume that
\begin{itemize}
\item 
we  know the  weight multiplicities  of  $\Galg$ for  all $L(\lambda)$  with
$p$-restricted $\lambda$,
\item
we can  compute representatives of  the conjugacy  classes of $G(q)$  in the
given representation,
\item
we can identify  the composition factors of this representation  in terms of
their labels by highest weights,
\item
and one of:
\begin{itemize}
\item
the ordinary  character table  of $G(q)$  is known and  we can  identify the
conjugacy classes  in the given  representation with those in  the character
table, or
\item 
we  can   compute  the  character  table   of  $G(q)$  (in  this   case  the
identification  of the  conjugacy classes  with those  of a  known table  is
automatic and essentially unique).
\end{itemize}
\end{itemize}

Sometimes we  may need a  tool, called the  {\MeatAxe}~\cite{MeatAxe}. Given
representing matrices of a set of generators of a group over a finite field,
it can  find representing  matrices of the  generators for  each (absolutely
irreducible) composition factor.

Given  representing  matrices   of  a  set  of  group   generators  for  two
representations  over the  same field  it  is easy  to compute  representing
matrices  for  the  tensor  product  representation  (Kronecker  product  of
matrices). If we want to compute  Brauer character values for group elements
in several  representations we  just add representatives  of all  classes of
$p'$-elements to our set of group generators.

The table of Brauer characters is now computed as follows.
\begin{itemize}
\item[(1)]
If the given representation of  $G(q)$ is not absolutely irreducible compute
the composition factors with the {\MeatAxe}.
\item[(2)]
Compute  the Brauer  character of  the  trivial representation,  and of  the
composition factors  we have found.  (Recall that we know  the corresponding
highest weights.)
\item[(3)]
Using weight multiplicities  we determine the composition  factors of tensor
products  of representations  for which  we already  know the  corresponding
Brauer characters, using the method sketched in~\ref{ss_computechars}.
\item[(3a)] 
If there is a tensor product  which contains only one composition factor for
which we  do not  yet know  the Brauer  character, then  we can  compute the
Brauer character  of this  composition factor  by tensoring  and subtracting
known Brauer characters.
\item[(3b)]
Otherwise, we  use a  tensor product  which has  one or  several composition
factors  which  are  easy  to  label  (e.g.,  via  their  degrees),  compute
representing  matrices  for  this  tensor product,  use  the  {\MeatAxe}  to
find  the composition  factors, and  compute  the Brauer  characters of  new
composition factors as above.
\item[(4)] 
If not all Brauer characters are found go back to step~(3).
\end{itemize}

We do not call  this description an algorithm, because it  is not clear that
we will  always be  able to  identify the label  of new  composition factors
found with the {\MeatAxe}.  And there can be a practical  problem if we need
to  apply the  {\MeatAxe} to  representations of dimensions exceeding a  few
thousands.

But in practice this method worked very well in cases we have tried.

\subsection{Example $G(q) = \Spin_8^-(3)$}\label{ss-2D4-method1}

\sloppy
In this  case we can find  a $16$-dimensional representation of  $G(q)$ over
$\F_3$ and we  can compute class representatives and  the ordinary character
table with {\Magma}~\cite{Magma}.

\fussy
The {\MeatAxe} yields two $8$-dimensional absolutely irreducible composition
factors over $\F_9$.

We have a description  of $G(q)$ in terms of a  root datum as in~\ref{ex-d4}
and  we know  the weight  multiplicities of  $\Galg$ for  all $3$-restricted
weights as stated in~\ref{ThmMults}.

Using  weight  multiplicities we  compute  (writing  $\chi_\lambda$ for  the
character of $L(\lambda)$):
\[
\chi_{(1,0,0,0)} \otimes \chi_{(0,1,0,0)} = \chi_{(1,1,0,0)} +
\chi_{(0,0,0,1)} 
\]
\[
\chi_{(1,0,0,0)} \otimes \chi_{(1,0,0,0)} = \chi_{(2,0,0,0)} +
\chi_{(0,0,1,0)} + \chi_{(0,0,0,0)}
\]

$\Galg$ has three  $8$-dimensional irreducible representations corresponding
to the  highest weights  $\lambda \in \{(1,0,0,0),  (0,1,0,0), (0,0,0,1)\}$.
Using  the  {\MeatAxe} we  find  the  third  irreducible  of degree  $8$  as
composition factor of the tensor product of  the given ones. And we find the
$28$-dimensional  module  $L((0,0,1,0))$  from  the  tensor  product  of  an
$8$-dimensional module with itself.

The  three  $8$-dimensional  representations   are  permuted  by  the  first
Frobenius  twist  as   described  by  $F_0$  (see~\ref{ss_computechars}(g)).
Applying  the  corresponding Galois  automorphism  to  the Brauer  character
values, we  find that the two  representations we started with  are swapped.
So, we started with $L((1,0,0,0))$ and $L((0,1,0,0))$ (in any order) and the
third $8$-dimensional module is $L((0,0,0,1))$.

At  this  stage we  know  the  Brauer  characters  of $L(\lambda)$  for  all
fundamental weights  $\lambda$. It turns  out that  from here we  can always
find a  tensor product of  known Brauer  characters which contains  only one
constituent whose Brauer character is not yet known:

\[\chi_{(0,0,0,2)} = \chi_{(0,0,0,1)} \otimes \chi_{(0,0,0,1)} -
\chi_{(0,0,1,0)} - \chi_{(0,0,0,0)}\]
\[\chi_{(0,1,0,1)} = \chi_{(0,0,0,1)} \otimes \chi_{(0,1,0,0)} -
\chi_{(1,0,0,0)} \]
\[ \ldots \]
\[\begin{array}{rcl}
\chi_{(2,1,0,2)} &=&\chi_{(1,1,0,2)} \otimes \chi_{(1,0,0,0)} -
\chi_{(0,1,1,2)} - 2\chi_{(1,2,0,1)} -2\chi_{(1,0,0,3)} \\
&& -3\chi_{(1,0,1,1)} -2 \chi_{(2,1,0,0)} - 4\chi_{(0,1,0,2)}
-2\chi_{(0,1,0,0)} \end{array}
\]
and so  on. That  is, we can  find all remaining  Brauer characters  just by
simple computations with characters.

\section{Semisimple Classes}\label{sec-sscl}

In  this section  we find  representatives of  semisimple conjugacy  classes
without an explicit  representation of $G(q)$. Let  $\Talg \subseteq \Galg$,
$F$, $F_0$, $F(\Talg) = \Talg$, $G(q)$ as in previous sections.

We use the following facts, see~\cite[3.7, 3.1, 3.2, 3.5]{Ca85}.

Each semisimple conjugacy of $\Galg$ intersects $\Talg$ and the intersection
is a single $W$-orbit of $\Talg$.

A semisimple  conjugacy class of  $\Galg$ intersects  $G(q)$ if and  only if
the  class  is  $F$-stable.  In  that case  the  intersection  is  a  single
$G(q)$-conjugacy class.

The maximal  torus can be  recovered from the root  datum as $\Talg  \cong Y
\otimes_\Z \Fp^\times$.

Writing $\Qpp$ for  the additive group of rational  numbers with denominator
not  divisible  by  $p$  there   is  an  isomorphism  $\overFp^\times  \cong
(\Qpp/\Z)^+$. There  is also an  explicit choice  of such an  isomorphism in
terms of  Conway polynomials,  such that  the composition    of $(\Qpp/\Z)^+
\to  \C^\times$,  $\frac{r}{s}+\Z  \mapsto  \exp(2  \pi  i  r/  s)$,  yields
the embedding  $\overFp^\times \hookrightarrow  \C^\times$ mentioned  in the
beginning of Section~\ref{sec-smallrep}.

The centralizer in $\Galg$  of an element $t \in \Talg$  is also a reductive
group, parametrized by the same lattices  $X$, $Y$ as $\Galg$ and the subset
of the  roots $\alpha  \in \Phi$ with  $\alpha(t) = 0  \in \Qpp/\Z$  and the
corresponding coroots.  If $w \in W$  with $w(F(t)) = t$  then the Frobenius
morphism on the centralizer is described by the matrix $F_0w$.

Combining  this and  using  the chosen  $\Z$-basis of  $Y$  we can  identify
$\Talg$ with  the additive group $(\Qpp/\Z)^l$.  The action of $w  \in W$ on
$Y$  extends to  an action  on $(\Qpp/\Z)^l$,  $w(t) =  t w$  (formal matrix
multiplication). The same holds for the Frobenius action, $F(t) = t (qF_0)$.
Evaluating a weight $\lambda  \in X$ on $t \in \Talg  = (\Qpp/\Z)^l$ is also
done by a matrix product $\lambda(t) = \lambda t^{tr}$.

This leads to the following  algorithm to determine the semisimple conjugacy
classes of $G(q)$ by finding representatives of the $F$-stable $W$-orbits of
$\Talg$.

\begin{itemize}
\item[(1)] 
Determine a set of representatives $F_0w \in F_0W$ under conjugation of $W$.
(Or, representatives $w$ of the $F_0$-conjugacy classes of $W$.)
\item[(2)]
For each $F_0w$ found in~(1) find all  solutions of the equation $t(q F_0w -
\id) = 0 \in (\Qpp/\Z)^l$.
\item[(3)]
For each element $t$ found in~(2) compute its $W$-orbit in $(\Qpp/\Z)^l$ and
take the (lexicographically) minimal element as representative.
\item[(4)]
For each  representative $t$  from~(3) compute the  roots $\alpha  \in \Phi$
with $\alpha(t) = 0 \in (\Qpp/\Z)$, and a $w \in W$ with $w(F(t)) = t$.
\end{itemize}

Note  that  we  may  find  some  orbits/classes  several  times  during  the
algorithm.

We  remark  that  in  practice  we  have  used  programs  which  parametrize
semisimple  classes generically  for $G(q)$  for  all prime  powers $q$  and
specialized to  the specific small  $q$ considered  here. But the  much more
elementary approach sketched above is sufficient for the application in this
paper.

In  addition  to the  representatives  of  semisimple  classes we  need  the
following information:
\begin{description}
\item[Power maps.]
For small positive integers $k$ we want to know for each semisimple class of
an  element $s  \in  G(q)$  the class  of  $s^k$. If  the  class  of $s$  is
represented by an element $t \in  (\Qpp/\Z)^l$ as above, we compute $kt$ and
the lexicographically minimal element of its $W$-orbit.
\item[Multiplication with central elements.] 
For each semisimple class of an element $s \in G(q)$ and each element $z \in
Z(G(q))$ we want  to know the class  of $sz$ (this is well  defined). If the
class is represented by $t \in (\Qpp/\Z)^l$ and $z$ is represented by $c \in
(\Qpp/\Z)^l$  we  compute  the  lexicographically  minimal  element  in  the
$W$-orbit of $t+c$.
\end{description}

\subsection{Example}\label{ex-sstorus}

We   use  again   the  example   $G(q)   =  \Spin_8^-(3)$   and  the   setup
from~\ref{ex-d4}.  The matrix  of $F_0$  is  the permutation  matrix of  the
transposition $(1,2)$. The computations can be done with the basic functions
of the {\CHEVIE} package~\cite{CHEVIE}.

Let $w  = s_{4}^\vee s_{3}^\vee s_{2}^\vee  s_{1}^\vee s_{3}^\vee s_{4}^\vee
s_{ 1}^\vee s_{3}^\vee s_{1}^\vee =
{\footnotesize \left(\begin{array}{rrrr}%
0&0&-1&0\\%
1&1&1&0\\%
0&-1&0&0\\%
-1&0&-1&-1\\%
\end{array}\right)}$ 
as matrix acting on $Y$. We have to consider the equation
\[t (qF_0^{tr} w -\id) = t M =0 \in (\Qpp/\Z)^l, \textrm{ where }
M = {\footnotesize \left(\begin{array}{rrrr}%
2&3&3&0\\%
0&-1&-3&0\\%
0&-3&-1&0\\%
-3&0&-3&-4\\%
\end{array}\right)}.\]
With the  Smith normal form algorithm  we find unimodular matrices  $L,R \in
\Z^{4 \times 4}$  with $ L M R =  \textrm{diag}(0,0,8,8)$. So, the solutions
of the original equation are
\[(0,0,i/8,j/8)\cdot  L \textrm{ with } 0 \leq i,j < 8.\]

Considering the  specific solution $t' =  (0,0,1/4,1/4) L = (1/4,  1/4, 1/2,
1/2)$ its $W$-orbit contains $8$  torus elements, the minimal representative
in that orbit is $t =  (1/4,1/4,0,0)$. The roots $\{\alpha \in \Phi\mid\; t'
\alpha = 0\}$ form a subsystem of  the root system of $\Galg$ of type $A_3$.
This  yields the  root  datum  of the  centralizer  $C  = C_\Galg(t')$.  The
Frobenius action on this centralizer  is described by the matrix $F_0w^{tr}$
(the transposed of  $w$ described the action  of the same element  of $W$ on
$X$). With {\CHEVIE}  we see that the  centralizer $C$ is of  type $A_3(q) +
T(q+1)$ (Dynkin  diagram of type $A_3$  with trivial Frobenius action  and a
central  torus  $Z^0$  with  $|(Z^0)^F|  = q+1$).  For  $q=3$  we  find  the
centralizer order $48522240$.

We see that $t'$ has order $4$ and  we can identify the classes of $kt'$ for
$k = 2,3$ and so all power maps for this class.

The center  of $G(q)$  is of order  $2$ and the  non-trivial element  in the
center  is $c  =  (1/2, 1/2,  0,  0)$. The  element $t'+c$  is  in the  same
$W$-orbit as $t'$.

\section{Brauer Table from Weight Multiplicities}\label{sec-brauer}
 
Let $\Galg$, $\Talg$, $p$, $F$, $G(q)$ as in the previous sections. We fix a
$q$ and assume that

\begin{itemize}
\item
we  know the  weight multiplicities  of  $\Galg$ for  all $L(\lambda)$  with
$p$-restricted $\lambda$,
\item
we know the ordinary character table  of $G(q)$ (abstractly, that is without
a  labelling of  the  conjugacy  classes by  representatives  in a  concrete
groups), including the power maps of the classes,
\item
we have the representatives of semisimple classes of $G(q)$ in form of torus
elements as explained in Section~\ref{sec-sscl}.
\end{itemize}

With this information it is easy  to compute the values of Brauer characters
as functions on  the given representatives of semisimple classes.  Let for a
fixed $p$-restricted dominant weight  $\lambda$ the weight multiplicities of
$L(\lambda)$ be given  as a dominant character, and let  $t \in (\Qpp/\Z)^l$
be a representative of a semisimple class.  If $\mu$ is a dominant weight of
$L(\lambda)$ with multiplicity $m_\mu$ then for each $\mu'$ in the $W$-orbit
of $\mu$ the element  $t$ has $m_\mu$ times the eigenvalue  $a = \mu' t^{tr}
\in  \Qpp/\Z$ on  the weight  space $L(\lambda)_{\mu'}$.  We can  lift these
eigenvalues to $\exp(2 \pi i a) \in  \C^\times$ and add them all up over all
weights  of $L(\lambda)$  to  find  the Brauer  character  value  of $t$  on
$L(\lambda)$.

To  be able  to relate  the  Brauer characters  found so  far with  ordinary
characters  we  need   to  identify  the  conjugacy   classes  described  by
representatives $t \in \Talg$ with the classes of $p'$-elements in the given
ordinary character table  of $G(q)$. This map is usually  not unique and can
be difficult to determine.

Recall  that our  class representatives  in  $\Talg$ found  as described  in
Section~\ref{sec-sscl} come with the following information:
\begin{itemize}
\item element order,
\item centralizer order,
\item power maps,
\item permutation of classes by multiplication with central elements.
\end{itemize}

All of this  information is also contained in the  abstract character tables
as they  come from the  {\ATLAS}~\cite{ATLAS} or the {\GAP}  character table
library {\CTblLib}~\cite{CTblLib} (the  permutation from multiplication with
central  elements  can  be  computed from  so  called  class  multiplication
coefficients).

Of  course,  the identification  of  classes  we  are  looking for  must  be
compatible with these data.

For  character tables  in {\GAP}  we can  compute the  group of  \emph{table
automorphisms}.  These consist  of  permutations of  the conjugacy  classes,
compatible with  power maps, which  leave the set of  irreducible characters
invariant. We want to find the identifications of classes modulo these table
automorphisms.

It is not a  priori clear that modulo table automorphisms  there is only one
such identification.  And we do  not have  a practical algorithm  which will
find all  possible identifcations  compatible with our  data (of  course, by
brute force one can try all  possibilities and check compatibility, but that
is not  practical). Therefore,  each case  needs some  ad hoc  procedure. We
mention some typical arguments in the example below.

\subsection{Example $G(q) = \Spin_8^-(3)$}\label{ss-ident-2d43}

We  consider  again   the  group  $G(q)  =  \Spin_8^-(3)$.   We  have  found
representatives $t \in \Talg$ of the $3^4 = 81$ semisimple classes of $G(q)$
together with  the data mentioned  above. The  character table of  $G(q)$ is
available in {\GAP} under the  name \texttt{"2.O8-(3)"}. The table has $640$
table automorphisms.

The identifications of the two center elements are clear.

There are for example $10$ classes of elements of order $82$ and the squares
of their representative  are in $10$ classes of elements  of order $41$ (the
$41$st powers are the non-trivial element  in the center). It turns out that
there is a table automorphism which permutes the $10$ classes of elements of
order $82$ cyclically and their $2$nd powers accordingly and fixes all other
classes. Furthermore, the $7$th power map  also permutes the $10$ classes of
$82$-elements cyclically.  This shows that we  can identify one class  of an
$82$-element arbitrarily and then the identification of all other classes of
$82$-elements and $41$-elements is determined from the $7$th and $2$nd power
maps.

A similar argument works for $8$ classes of $104$-elements and their powers.
The  choice of  the  identification of  one  of these  classes  can be  made
independently  from the  choice  for the  $82$-elements  because the  common
powers of these elements are only the center elements.

There  are   $4$  classes  of   $40$-elements  with  centralizer   of  order
$160$.  Their $5$th  powers  yield  $8$-elements which  are  also powers  of
$104$-elements  and  so  are  already   identified.  This  leaves  only  one
possibility for identifying these $4$ classes.

There are $3$ classes of $8$-elements  with centralizer order $64$. Only one
of  these classes  has  elements  whose $2$nd  power  has centralizer  order
$69120$ which  fixes the  identification of  that class.  For the  other two
classes there  are two  possible identifications  which are  both compatible
with our data.

We can  proceed with  this type  of arguments until  the number  of possible
identifications  becomes reasonably  small:  two choices  for  one class  of
$8$-elements as mentioned, four choices  each for one class of $40$-elements
and for  one class of  $56$-elements, and two  choices for another  class of
$40$-elements.

At this stage  we just try out  any of the $64$ combinations  of choices. In
each case  the full identification  of classes  then follows from  the power
maps.  For each  of  these indentifications  we  construct the  hypothetical
Brauer character  table and  compute the corresponding  decomposition matrix
(that is express the restrictions of the ordinary characters to $p'$-classes
as  linear  combinations  of  the  Brauer  characters).  The  entries  of  a
decomposition matrix  must be  non-negative integers. This  conditions rules
out $60$  of the  $64$ possibilities  because these  yield some  non integer
or  negative  coefficients.  Two  pairs   of  the  remaining  four  possible
identifications differ  only by a  table automorphism, so that  modulo table
automorphisms we are left with two possible identifications.

Further  investigation  shows that  from  one  of  the two  possible  Brauer
character tables  we get the  other via  a Galois automorphism  which raises
complex roots of  unity to their $647$th power. This  shows that both tables
are  correct with  respect to  some choice  of $p$-modular  system, or  with
respect to some choice of the identification $\overFp^\times \cong \Qpp/\Z$.

\subsection{Some new Brauer character tables}\label{ss-new-brauer}

Using  the  techniques  described  in   this  paper  we  get  the  following
contribution to the Modular Atlas Project.

\begin{Thm}
The  $p$-modular  Brauer   character  tables  and  their   fusion  into  the
corresponding ordinary character table are known for the following cases:
{\rm
\[ \begin{array}{lllll}
p& \textrm{Lie type} & \textrm{group name} & 
\textrm{{\ATLAS} name} & \textrm{name in {\GAP}}\\
\hline
2 & F_4 & F_4(2) & F_4(2) &\texttt{F4(2)}\\
3 & D_4 & \Spin^+_8(3) & 2^2.O_8^+(3) & \texttt{2\^{}2.O8+(3)}\\
3 & D_4 & \Spin^-_8(3) & 2.O_8^-(3) & \texttt{2.O8-(3)}\\
2 & D_5 & \Spin^+_{10}(2) & O_{10}^+(2) & \texttt{O10+(2)}\\
2 & D_5 & \Spin^-_{10}(2) & O_{10}^-(2) & \texttt{O10-(2)}\\
2 & C_5 & \Sp_{10}(2) & S_{10}(2) & \texttt{S10(2)}\\
\end{array}
\]
}

Furthermore,  partial   $2$-modular  Brauer   character  tables   are  known
for   the   groups   $E_6(2)$    and   $^2\!E_6(2)_{sc}$   ({\ATLAS}   names
$E_6(2)$   and  $3.{}^2\!E_6(2)$,   {\GAP}  names   {\rm\texttt{E6(2)}}  and
{\rm\texttt{3.2E6(2)}}). In these cases $44$  of the $64$ irreducible Brauer
characters are known. (We cannot get decomposition numbers from this partial
information.) \end{Thm}

These character tables are too big to  be printed in this article. They will
be  available in  future  versions  of the  {\GAP}  character table  library
{\CTblLib}~\cite{CTblLib}. 

Actually,  not all  ordinary character  tables mentioned  above are  printed
in  the  {\ATLAS}~\cite{ATLAS}  (but  they are  available  in  {\GAP}):  for
$2^2.O_8^+(3)$ only the  table of the simple quotient is  printed, tables of
$2.O_8^-(3)$ or its simple quotient are not printed, the availability of the
table of  $S_{10}(2)$ is only  mentioned in the \textit{Improvements  to the
ATLAS}~\cite[App.  2]{ModAtl}. For  $3.{}^2\!E_6(2)$ only  the table  of the
simple quotient was printed, the table  of the extension was computed by the
author~\cite{L2E6}. The table  of $E_6(2)$ was not printed  in the {\ATLAS}.
Actually,  while trying  to identify  semisimple  classes of  this table  we
discovered an  error in that ordinary  table. It turned out  that Bill Unger
had recently  recomputed that  table with {\Magma}  and also  discovered the
error. A  corrected ordinary  table  of $E_6(2)$ will also  become available
with future version of the character table library {\CTblLib}.

\bibliographystyle{alpha}
\bibliography{Luebeck_BrauerDefChar}

\end{document}